\documentclass[11pt,a4paper]{amsart}%
\pdfoutput=1
\usepackage{hyperref}
\usepackage{amsfonts}
\usepackage{amsmath}
\usepackage{amssymb}
\usepackage{graphicx}
\usepackage{algorithm}
\usepackage{algorithmic}%
\newtheorem{theorem}{Theorem}

\theoremstyle{definition}
\newtheorem{definition}[theorem]{Definition}
\newtheorem{proposition}[theorem]{Proposition}
\newtheorem{example}[theorem]{Example}

\theoremstyle{remark}

\numberwithin{equation}{section}
\theoremstyle{plain}

\begin{document}
\title[The Kustin-Miller complex construction]{Implementing the Kustin-Miller complex construction}
\author{Janko B\"{o}hm}
\address{Department of Mathematics, Universit\"{a}t des Saarlandes, Campus E2 4 \\
D-66123 \\
Saarbr\"{u}cken, Germany}
\email{boehm@math.uni-sb.de}
\author{Stavros Argyrios Papadakis}
\address{Centro de An\'{a}lise Matem\'{a}tica, Geometria e Sistemas Din\^{a}micos,
Departamento de Matem\'atica, Instituto Superior T\'ecnico, Universidade
T\'ecnica de Lisboa, Av. Rovisco Pais, 1049-001 Lisboa, Portugal}
\email{papadak@math.ist.utl.pt}
\thanks{J. B. was supported by DFG (German Research Foundation) through Grant
BO3330/1-1. S. P. was supported by the Portuguese Funda\c{c}\~ao para a
Ci\^encia e a Tecno\-lo\-gia through Grant SFRH/BPD/22846/2005 of
POCI2010/FEDER and through Project PTDC/MAT/099275/2008. }

\begin{abstract}
The Kustin-Miller complex construction, due to A. Kustin and M. Miller, can be
applied to a pair of resolutions of Gorenstein rings with certain properties
to obtain a new Gorenstein ring and a resolution of it. It gives a tool to
construct and analyze Gorenstein rings of high codimension. We describe the
Kustin-Miller complex and its implementation in the Macaulay2 package
\textsc{KustinMiller}, and explain how it can be applied to explicit examples.

\end{abstract}
\maketitle

\section{Introduction\label{Sec Introduction}}

Many important rings in commutative algebra and algebraic geometry turn out to
be Gorenstein rings, i.e., commutative rings such that the localization at
each prime ideal is a Noetherian local ring $R$ with finite injective
dimension as an $R$-module. Examples are canonical rings of regular algebraic
surfaces of general type, anticanonical rings of Fano varieties and
Stanley-Reisner rings of triangulations of spheres. Except for the complete
intersection cases of codimension $1$ and $2$ a structure theorem for
Gorenstein rings is known only for codimension $3$ by the theorem of
Buchsbaum-Eisenbud \cite{BE}, which describes them in terms of Pfaffians of a
skew-symmetric matrix. One goal of unprojection theory, which was introduced
by A. Kustin, M. Miller and M. Reid and developed further by the second author
(see, e.g., \cite{KM}, \cite{R1}, \cite{PR}, \cite{P1}), is to act as a
substitute for a structure theorem in codimension $\geq4$ by providing a
construction to increase the codimension in a non-trivial way, while staying
in the class of Gorenstein rings. The geometric motivation is to provide
inverses of certain projections in birational geometry. The process can be
considered as a version of Castelnuovo blow-down.

Examples of applications range from the construction of Campedelli surfaces
\cite{NP} to results on the structure of Stanley-Reisner rings \cite{BPcyclic}%
. For an outline of more applications see \cite{R1}, the introduction of
\cite{BP} and Section \ref{Sec Applications} below.

We describe the Kustin-Miller complex construction \cite{KM}, which is the key
tool to obtain resolutions of unprojection rings, and discuss our
implementation in the \textsc{Macaulay2} \cite{GS} package
\textsc{KustinMiller} \cite{BP2}. We illustrate the construction with examples
and applications.

\section{Implementation of the Kustin-Miller complex
construction\label{Sec Kustin Miller}}

We will consider the following setup: Let $R$ be a positively graded
polynomial ring over a field and $I,J\subset R$ homogeneous ideals of $R$ such
that $R/I$ and $R/J$ are Gorenstein, $I\subset J$ and $\dim R/J=\dim R/I-1$.
By \cite[Proposition~3.6.11]{BH} there are $k_{1},k_{2}\in\mathbb{Z}$ such
that $\omega_{R/I}=R/I(k_{1})$ and $\omega_{R/J}=R/J(k_{2})$. Assume that
$k_{1}>k_{2}$ so that the unprojection ring defined below is also positively graded.

\begin{definition}
\label{def unproj ring}\cite{PR} Let $\phi\in\operatorname{Hom}_{R/I}(J,R/I)$
be a homomorphism of degree $k_{1}-k_{2}$ such that $\operatorname{Hom}%
_{R/I}(J,R/I)$ is generated as an $R/I$-module by $\phi$ and the inclusion
morphism $i$. We call the graded algebra $R[T]/U$, where $T$ is a variable of
degree $k_{1}-k_{2}$ and%
\[
U=\left(  I,\;Tu-\phi(u)\mid u\in J\right)
\]
the \textbf{Kustin--Miller unprojection ring} of the pair $I\subset J$ defined
by $\phi$.
\end{definition}

\begin{proposition}
\cite{KM,PR} The ring $R[T]/U$ is Gorenstein and independent of the choice of
$\phi$ (up to isomorphism).
\end{proposition}

Following \cite{KM}, we now describe the construction of a graded free
resolution of $R[T]/U$ from those of $R/I$ and $R/J$. We will refer to this as
the \textbf{Kustin-Miller complex construction}. Denote by $g=\dim R-\dim R/J$
the codimension of the ideal $J$ of $R$, and suppose $g\geq4$ (the special
cases $g=2$ and $3$ can be treated in a similar way). Let
\[%
\begin{tabular}
[c]{ll}%
$C_{J}:$ & $R/J\leftarrow A_{0}\overset{a_{1}}{\leftarrow}A_{1}\overset{a_{2}%
}{\leftarrow}\dots\overset{a_{g-1}}{\leftarrow}A_{g-1}\overset{a_{g}%
}{\leftarrow}A_{g}\leftarrow0$\\
$C_{I}:$ & $R/I\leftarrow B_{0}\overset{b_{1}}{\leftarrow}B_{1}\overset{b_{2}%
}{\leftarrow}\dots\overset{b_{g-1}}{\leftarrow}B_{g-1}\leftarrow0$%
\end{tabular}
\
\]
be minimal graded free resolutions (self-dual by the Gorenstein property,
\cite{Ei}) of $R/J$ and $R/I$ as $R$-modules with $A_{0}=B_{0}=R$,
$A_{g}=R\left(  k_{1}-\eta\right)  $ and $B_{g-1}=R\left(  k_{2}-\eta\right)
$, where $\eta$ is the sum of the degrees of the variables of $R$. Consider
the complex%
\[%
\begin{tabular}
[c]{ll}%
$C_{U}:$ & $R[T]/U\leftarrow F_{0}\overset{f_{1}}{\leftarrow}F_{1}%
\overset{f_{2}}{\leftarrow}\dots\overset{f_{g-1}}{\leftarrow}F_{g-1}%
\overset{f_{g}}{\leftarrow}F_{g}\leftarrow0$%
\end{tabular}
\
\]
with the modules%
\begin{gather*}%
\begin{tabular}
[c]{lll}%
$F_{0}=B_{0}^{\prime}$, &  & $F_{1}=B_{1}^{\prime}\oplus A_{1}^{\prime}%
(k_{2}-k_{1})$%
\end{tabular}
\\%
\begin{tabular}
[c]{ll}%
$F_{i}=B_{i}^{\prime}\oplus A_{i}^{\prime}(k_{2}-k_{1})\oplus B_{i-1}^{\prime
}(k_{2}-k_{1})$ & f$\text{or }2\leq i\leq g-2$%
\end{tabular}
\\%
\begin{tabular}
[c]{lll}%
$F_{g-1}=A_{g-1}^{\prime}(k_{2}-k_{1})\oplus B_{g-2}^{\prime}(k_{2}-k_{1})$, &
& $F_{g}=B_{g-1}^{\prime}(k_{2}-k_{1})$%
\end{tabular}
\end{gather*}
where for an $R$-module $M$ we denote $M^{\prime}:=M\otimes_{R}R[T]$.

By specifying chain maps $\alpha:C_{I}\rightarrow C_{J}$, $\beta
:C_{J}\rightarrow C_{I}[-1]$ and a homotopy map (not necessarily chain map)
$h:C_{I}\rightarrow C_{I}$ we will define the differentials as%
\begin{gather*}%
\begin{tabular}
[c]{lll}%
$f_{1}=\left(
\begin{array}
[c]{cc}%
b_{1} & \beta_{1}+T\cdot a_{1}%
\end{array}
\right)  $, &  & $f_{2}=\left(
\begin{array}
[c]{ccc}%
b_{2} & \beta_{2} & h_{1}+T\cdot I_{1}\\
0 & -a_{2} & -\alpha_{1}%
\end{array}
\right)  $%
\end{tabular}
\\%
\begin{tabular}
[c]{ll}%
$f_{i}=\left(
\begin{array}
[c]{ccc}%
b_{i} & \beta_{i} & h_{i-1}+(-1)^{i}T\cdot I_{i-1}\\
0 & -a_{i} & -\alpha_{i-1}\\
0 & 0 & b_{i-1}%
\end{array}
\right)  $ & $\text{for }3\leq i\leq g-2$%
\end{tabular}
\\%
\begin{tabular}
[c]{ll}%
$\hspace{-0.05in}f_{g-1}=\left(
\begin{array}
[c]{cc}%
\beta_{g-1} & h_{g-2}+(-1)^{g-1}T\cdot I_{g-2}\\
-a_{g-1} & -\alpha_{g-2}\\
0 & b_{g-2}%
\end{array}
\right)  $, & $f_{g}=\left(
\begin{array}
[c]{c}%
-\alpha_{g-1}+(-1)^{g}\frac{1}{\beta_{g}\left(  1\right)  }T\cdot a_{g}\\
b_{g-1}%
\end{array}
\right)  $%
\end{tabular}
\end{gather*}
where $I_{t}$ denotes the $\operatorname{rank}B_{t}\times\operatorname{rank}%
B_{t}$ identity matrix. We now discuss the construction of $\alpha$, $\beta$
and $h$:

Fix $R$-module bases $e_{1},\dots,e_{t_{1}}$ of $A_{1}$ and $\hat{e}_{1}%
,\dots,\hat{e}_{t_{1}}$ of $A_{g-1}$ and write%
\[%
\begin{tabular}
[c]{lll}%
$\sum_{i=1}^{t_{1}}\hat{c}_{i}\cdot\hat{e}_{i}:=a_{g}(1_{R})$, &  &
$c_{i}\cdot1_{R}:=a_{1}(e_{i})\text{ for }i=1,...,t_{1}$%
\end{tabular}
\
\]
where by Gorensteinness $c_{i},\hat{c}_{i}\in J$ for all $i$. Denote by
$l_{i},\hat{l}_{i}\in R$ lifts of $\phi(c_{i}),\phi(\hat{c}_{i})\in R/I$,
respectively. For an $R$-module $A$ we write $A^{\ast}=\operatorname{Hom}%
_{R}(A,R)$ and for an $R$-basis $f_{1},\dots f_{t}$ of $A$ we denote by
$f_{1}^{\ast},\dots,f_{t}^{\ast}$ the dual basis of $A^{\ast}$. Now consider
the $R$-homomorphism
\[
A_{g-1}^{\ast}\rightarrow R=B_{g-1}^{\ast}\text{,\quad}\hat{e}_{i}^{\ast
}\mapsto\hat{l}_{i}\cdot1_{R}%
\]
which (by self-duality of $C_{I},C_{J}$) extends to a chain map $C_{J}^{\ast
}\rightarrow C_{I}^{\ast}$ and denote by $\tilde{\alpha}:C_{I}\rightarrow
C_{J}$ its dual. The map $\tilde{\alpha}_{0}:B_{0}=R\rightarrow R=A_{0}$ is
multiplication by an invertible element of $R$, cf. \cite{P1}, set
$\alpha=\tilde{\alpha}/\tilde{\alpha}_{0}\left(  1_{R}\right)  $.

We obtain $\beta:C_{J}\rightarrow C_{I}[-1]$ by extending%
\[
\beta_{1}:A_{1}\rightarrow R=B_{0}\text{,\quad}e_{i}\mapsto-l_{i}\cdot1_{R}%
\]

Finally, by \cite[p.~308]{KM} there is a homotopy $h:C_{I}\rightarrow C_{I}$
with $h_{0}=h_{g-1}=0$ and
\[
\beta_{i}\alpha_{i}=h_{i-1}b_{i}+b_{i}h_{i}\text{\quad for }1\leq i\leq g
\]

\begin{theorem}
\cite{KM} The complex $C_{U}$ is a graded free resolution of $R[T]/U$ as an
$R[T]$-module.
\end{theorem}

It is important to remark that $C_{U}$ is not necessarily minimal, although in
many examples coming from algebraic geometry it is.

We now describe, how to compute $C_{U}$, as implemented in our
\textsc{Macaulay2} package \textsc{KustinMiller}. First note, that we can
determine $\phi$ via the commands \texttt{Hom(J,R}\symbol{94}\texttt{1/I)} and
\texttt{homomorphism} available in \textsc{Macaulay2} . Furthermore one can
extend homomorphisms to chain maps by the command \texttt{extend}.

\begin{algorithm}                      
\caption{Kustin-Miller complex}          
\label{alg km}
\begin{algorithmic}[1]
\REQUIRE Resolutions $C_{I}$ and $C_{J}$, denoted as above, for homogeneous
ideals $I \subset J$ in a polynomial ring $R$ with $R/I$ and $R/J$ Gorenstein,
and $\dim R/J=\dim R/I-1$.
\ENSURE The Kustin-Miller complex $C_{U}$ associated to $I$ and $J$.
\STATE Compute $\phi$ as in Definition \ref{def unproj ring} above.
\STATE Compute the dual $C_{J}^{\ast}$ of $C_{J}$ and express the first
differential as the product of a square matrix $Q$ with $a_{1}$. Extend the
homomorphism $\phi\circ Q$ to a chain map $\alpha^{\ast}:C_{J}^{\ast
}\rightarrow C_{I}^{\ast}$ and dualize to obtain $\tilde{\alpha}%
:C_{I}\rightarrow C_{J}$. Dividing all differentials of $\tilde{\alpha}$ by
the inverse of the entry of $\tilde{\alpha}_{0}$ yields $\alpha:C_{I}%
\rightarrow C_{J}$.
\STATE Extend the map $A_{1}\rightarrow B_{0}$ given by $\phi$ to a chain map
$C_{J}\rightarrow C_{I}[-1]$ and multiply the differentials by $-1$ to obtain
$\beta:C_{J}\rightarrow C_{I}[-1]$.
\STATE Set $h_{0}:=0_{R}$.
\FOR{$i=1$ to $g-1$}
\STATE Set $h_{i}^{\prime}:=\beta_{i}\alpha_{i}-h_{i-1}b_{i}$.
\STATE Using the \texttt{extend} command obtain $h_{i}$ in the diagram%
\[%
\begin{tabular}
[c]{rll}%
$B_{i}$ & $\overset{h_{i}^{\prime}}{\longrightarrow}$ & $B_{i}$\\
${\small id}\uparrow\,$ &  & $\uparrow{\small b}_{i}$\\
$B_{i}$ & $\underset{h_{i}}{\longrightarrow}$ & $B_{i}$%
\end{tabular}
\]
\ENDFOR
\RETURN the differentials $f_{i}$ according to the formulas given above.
\end{algorithmic}
\end{algorithm}\vspace{-0.05in}

\section{Applications\label{Sec Applications}}

We comment on some applications of the Kustin-Miller complex construction
involving the authors (for examples on these, see the documentation of our
\textsc{Macaulay2} package \textsc{KustinMiller}).

\subsection{Cyclic polytopes}

For a polynomial ring $R=k\left[  x_{1},...,x_{n}\right]  $ denote by
$I_{d}\left(  R\right)  $ the Stanley-Reisner ideal of the boundary complex of
the cyclic polytope of dimension $d$ with vertices $x_{1},...,x_{n}$. As shown
in \cite{BPcyclic} the Kustin-Miller complex construction yields a recursion
for a minimal resolution of $I_{d}\left(  R\right)  $: For $d$ even apply
Algorithm \ref{alg km} with $T=x_{n}$ to minimal resolutions $C_{I}$ and
$C_{J}$ of $I=I_{d}\left(  k\left[  x_{1},...,x_{n-1}\right]  \right)  $ and
$J=I_{d-2}\left(  k\left[  z,x_{2},...,x_{n-2}\right]  \right)  $ considered
as ideals in $k\left[  z,x_{1},...,x_{n-1}\right]  $ and quotient by $\left(
z\right)  $. For $d$ odd one can proceed in a similar way.

\subsection{Stellar subdivisions}

Suppose $C$ is a Gorenstein* simplicial complex on the variables of $k\left[
x_{1},...,x_{n}\right]  $ and $F$ is a face of $C$. Let $C_{F}$ be obtained by
the stellar subdivision of $C$ with respect to $F$, introducing the new
variable $x_{n+1}$. Denote by $I$ the image of the Stanley-Reisner ideal of
$C$ in $k\left[  z,x_{1},...,x_{n}\right]  $ and by $J=\left(  z\right)  +I:(%
{\textstyle\prod\nolimits_{i\in F}}
x_{i})$ the ideal corresponding to the link of $F$. Apply Algorithm
\ref{alg km} to minimal resolutions of $I$ and $J$ with $T=x_{n+1}$ and
quotient by $\left(  z\right)  $. By \cite{BP} this yields a resolution of the
Stanley-Reisner ring of $C_{F}$.

\subsection{Constructions in Algebraic Geometry}

In the paper \cite{NP} a series of Kustin-Miller unprojections was used in
order to give the first examples of Campedelli algebraic surfaces of general
type with algebraic fundamental group $\mathbb{Z}/6$, while a similar
technique produced in \cite{NP2} seven families of Calabi-Yau $3$-folds of
high codimension. In both cases, the Kustin--Miller complex construction was
used to control the numerical invariants of the new varieties.

\section{Example\label{sec examples}}

\begin{example}
Using our \textsc{Macaulay2} package \textsc{KustinMiller} \cite{BP2} we
discuss an example given in \cite{P1} passing from a codimension $3$ to a
codimension $4$ ideal. Over the polynomial ring

\noindent\hspace{-0.08in}%
\begin{tabular}
[c]{cr}%
\texttt{i1:}\hspace{-0.08in} & \texttt{R = QQ[x\_1..x\_4, z\_1..z\_4];}%
\end{tabular}

\noindent consider the skew-symmetric matrix

\noindent\hspace{-0.08in}%
\begin{tabular}
[c]{crrrrrrrl}%
\texttt{i2:}\hspace{-0.08in} & \texttt{b2 = matrix\{} & \texttt{\{} &
\texttt{0,} & \texttt{x\_1,} & \texttt{x\_2,} & \texttt{x\_3,} & \texttt{x\_4}
& \texttt{\},}\\
&  & \texttt{\{} & \texttt{-x\_1,} & \texttt{0,} & \texttt{0,} &
\texttt{z\_1,} & \texttt{z\_2} & \texttt{\},}\\
&  & \texttt{\{} & \texttt{-x\_2,} & \texttt{0,} & \texttt{0,} &
\texttt{z\_3,} & \texttt{z\_4} & \texttt{\},}\\
&  & \texttt{\{} & \texttt{-x\_3,} & \texttt{-z\_1,} & \texttt{-z\_3,} &
\texttt{0,} & \texttt{0} & \texttt{\},}\\
&  & \texttt{\{} & \texttt{-x\_4,} & \texttt{-z\_2,} & \texttt{-z\_4,} &
\texttt{0,} & \texttt{0} & \texttt{\} \};}%
\end{tabular}

\noindent The Buchsbaum-Eisenbud complex

\noindent\hspace{-0.08in}%
\begin{tabular}
[c]{cc}%
\texttt{i3:}\hspace{-0.08in} & \texttt{betti( cI = resBE b2)}%
\end{tabular}

\noindent\hspace{-0.08in}%
\begin{tabular}
[c]{crcccc}
&  & 0 & 1 & 2 & 3\\
\texttt{o3:}\hspace{-0.08in} & \texttt{total:} & 1 & 5 & 5 & 1\\
& 0: & 1 & . & . & .\\
& 1: & . & 5 & 5 & .\\
& 2: & . & . & . & 1
\end{tabular}

\noindent resolves the ideal $I=\left(  b_{1}\right)  \subset R$ generated by
the $4\times4$-Pfaffians

\noindent\hspace{-0.08in}%
\begin{tabular}
[c]{cc}%
\texttt{i4:}\hspace{-0.08in} & \texttt{b1 = cI.dd\_1}%
\end{tabular}

\noindent\hspace{-0.08in}%
\begin{tabular}
[c]{cc}%
\texttt{o4:}\hspace{-0.08in} & \texttt{%
$\vert$%
z\_2z\_3-z\_1z\_4,-x\_4z\_3+x\_3z\_4,x\_4z\_1-x\_3z\_2,x\_2z\_2-x\_1z\_4,-x\_2z\_1+x\_1z\_3%
$\vert$%
}%
\end{tabular}

\noindent of the skew-symmetric matrix $b_{2}$. Consider the unprojection
locus $J$ with Koszul resolution

\noindent\hspace{-0.08in}%
\begin{tabular}
[c]{cc}%
\texttt{i5:}\hspace{-0.08in} & \texttt{J = ideal(z\_1..z\_4);}%
\end{tabular}

\noindent\hspace{-0.08in}%
\begin{tabular}
[c]{cc}%
\texttt{i6:}\hspace{-0.08in} & \texttt{betti( cJ = res J)}%
\end{tabular}

\noindent\hspace{-0.08in}%
\begin{tabular}
[c]{crccccc}
&  & \texttt{0} & \texttt{1} & \texttt{2} & \texttt{3} & \texttt{4}\\
\texttt{o6:}\hspace{-0.08in} & \texttt{total:} & \texttt{1} & \texttt{4} &
\texttt{6} & \texttt{4} & \texttt{1}\\
& \texttt{0:} & \texttt{1} & \texttt{4} & \texttt{6} & \texttt{4} & \texttt{1}%
\end{tabular}

\noindent Applying Algorithm \ref{alg km} we obtain the Kustin-Miller
resolution of the unprojection ideal $U\subset R\left[  T\right]  $, in this
case the ideal of the Segre embedding $\mathbb{P}^{2}\times\mathbb{P}%
^{2}\hookrightarrow\mathbb{P}^{8}$,

\noindent\hspace{-0.08in}%
\begin{tabular}
[c]{cc}%
\texttt{i7:}\hspace{-0.08in} & \texttt{betti( cU = kustinMillerComplex(cI, cJ,
QQ[T]))}%
\end{tabular}

\noindent\hspace{-0.08in}%
\begin{tabular}
[c]{crccccc}
&  & 0 & 1 & 2 & 3 & 4\\
\texttt{o7:}\hspace{-0.08in} & \texttt{total:} & 1 & 9 & 16 & 9 & 1\\
& 0: & 1 & . & . & . & .\\
& 1: & . & 9 & 16 & 9 & .\\
& 2: & . & . & . & . & 1
\end{tabular}

\noindent with generators

\noindent\hspace{-0.08in}%
\begin{tabular}
[c]{cc}%
\texttt{i8:}\hspace{-0.08in} & \texttt{f1 = cU.dd\_1}%
\end{tabular}%
\[
\left(  b_{1},\text{ }-x_{1}x_{3}+T\cdot z_{1},\text{ }-x_{1}x_{4}+T\cdot
z_{2},\text{ }-x_{2}x_{3}+T\cdot z_{3},\text{ }-x_{2}x_{4}+T\cdot
z_{4}\right)
\]
and syzygy matrix

\noindent\hspace{-0.08in}%
\begin{tabular}
[c]{cc}%
\texttt{i9:}\hspace{-0.08in} & \texttt{f2 = cU.dd\_2}%
\end{tabular}%
\[
\left(
\begin{tabular}
[c]{c|cccccc|ccccc}
& $0$ & $0$ & $0$ & $0$ & $0$ & $0$ & $T$ & $0$ & $0$ & $0$ & $0$\\
& $0$ & $0$ & $-x_{1}$ & $0$ & $0$ & $x_{2}$ & $0$ & $T$ & $0$ & $0$ & $0$\\
$b_{2}$ & $-x_{1}$ & $0$ & $0$ & $-x_{2}$ & $0$ & $0$ & $0$ & $0$ & $T$ & $0$
& $0$\\
& $0$ & $0$ & $-x_{3}$ & $-x_{3}$ & $-x_{4}$ & $0$ & $-x_{3}$ & $0$ & $0$ &
$T$ & $0$\\
& $0$ & $x_{3}$ & $0$ & $0$ & $0$ & $0$ & $0$ & $0$ & $0$ & $0$ & $T$\\\hline
& $z_{2}$ & $z_{3}$ & $0$ & $z_{4}$ & $0$ & $0$ & $z_{4}$ & $0$ & $-x_{4}$ &
$0$ & $x_{2}$\\
$0$ & $-z_{1}$ & $0$ & $z_{3}$ & $0$ & $z_{4}$ & $0$ & $0$ & $0$ & $x_{3}$ &
$-x_{2}$ & $0$\\
& $0$ & $-z_{1}$ & $-z_{2}$ & $0$ & $0$ & $z_{4}$ & $-z_{2}$ & $x_{4}$ & $0$ &
$0$ & $-x_{1}$\\
& $0$ & $0$ & $0$ & $-z_{1}$ & $-z_{2}$ & $-z_{3}$ & $0$ & $-x_{3}$ & $0$ &
$x_{1}$ & $0$%
\end{tabular}
\ \right)
\]

\end{example}

\noindent The code computing this example and various others related to the
applications mentioned above can be found in the documentation of the package
\textsc{KustinMiller} \cite{BP2}.

\end{document}